\documentclass[a4paper,12pt]{article}
\usepackage{comment}
\usepackage{cite}
\usepackage{amsmath}
\usepackage{amssymb}
\usepackage{amsfonts}
\usepackage[T1]{fontenc}
\usepackage[utf8]{inputenc}
\usepackage{graphicx}
\usepackage{fancyhdr}
\usepackage{float}
\usepackage{xcolor}
\usepackage{authblk}
\usepackage{mathrsfs}
\usepackage{empheq}
\usepackage[hyphens]{url}
\usepackage{hyperref} 
\usepackage[]{breakurl}
\usepackage{empheq}

\usepackage{geometry}
\begin{document}

\title{Series representations for the logarithm of the Glaisher-Kinkelin constant}

\author[$\dagger$]{Jean-Christophe {\sc Pain}$^{1,2,}$\footnote{jean-christophe.pain@cea.fr}\\
\small
$^1$CEA, DAM, DIF, F-91297 Arpajon, France\\
$^2$Universit\'e Paris-Saclay, CEA, Laboratoire Mati\`ere en Conditions Extr\^emes,\\ 
91680 Bruy\`eres-le-Ch\^atel, France
}

\maketitle

\begin{abstract}
In this note, we propose two series expansions of the logarithm of the Glaisher-Kinkelin constant. The relations are obtained using expressions of derivatives of the Riemann zeta function, and one of them involves hypergeometric functions.
\end{abstract}

\section{Introduction}

The Glaisher-Kinkelin constant $A$ is defined by \cite{Kinkelin1860,Glaisher1878,Glaisher1893,Glaisher1894,Voros1987}:
\begin{equation}
\lim_{n\rightarrow \infty}\frac{H(n)}{n^{n^2/2+n/2+1/12}}e^{n^2/4}=A, 	
\end{equation}
where $H(n)$ is the hyperfactorial \cite{Sloane1995}:
\begin{equation}
H(n)=\prod_{k=1}^nk^k.    
\end{equation}	
One has in particular
\begin{equation}\label{zeta1}
\ln A=\frac{1}{12}-\zeta'(-1)
\end{equation}
where $\zeta$ is the usual zeta function \cite{Edwards1974}, and
\begin{equation}\label{zeta2}
\sum_{k=2}^{\infty}\frac{\ln k}{k^2}=-\zeta'(2)={\frac{\pi^{2}}{6}}\left[12\ln A-\gamma-\ln(2\pi)\right],
\end{equation}
$\gamma$ being the Euler-Mascheroni constant. The Glaisher-Kinkelin constant is also given by \cite{Glaisher1878,Glaisher1894}:
\begin{equation}
A=2^{1/36}\pi^{1/6}e^{(-\gamma/4+s)/3}, 	
\end{equation}
where
\begin{equation}
s=\sum_{r=2}^{\infty}\frac{\left(-\frac{1}{2}\right)^r(2^r-1)\zeta(r)}{1+r}=\frac{1}{12}\left[3+3\gamma-36\zeta'(-1)-\ln 2-6 \ln \pi\right].	
\end{equation}
The calculation of the derivatives of the Riemann zeta function has been the subject of a number of studies (see for instance \cite{Apostol1985,Choudhury1995,Musser2011,Connon2015}). Apostol obtained the following expression, valid for $x>-2$ \cite{Apostol1985}:
\begin{equation}\label{derivzeta}
\zeta'(x)=-\frac{1}{(x-1)^2}+\frac{1}{12}-\frac{x(x+1)(x+2)}{6}I'_3(x)-\frac{(3x^2+6x+2)}{6}I_3(x),
\end{equation}
where
\begin{equation}\label{im}
I_m(s)=\int_1^{\infty}\frac{P_m(x)}{x^{s+m}}\,dx
\end{equation}
and
\begin{equation}\label{imderiv}
I'_m(s)=-\int_1^{\infty}\frac{P_m(x)}{x^{s+m}}\ln x\,dx.
\end{equation}
In section \ref{sec2}, we combine expression (\ref{zeta1}) with Apostol's formula (\ref{derivzeta}) for the derivative of the zeta function to obtain an expression of the Glaisher-Kinkelin constant. In section \ref{sec3}, we apply the same procedure to expression (\ref{zeta2}) to derive a second relation.

\section{Derivation of the first summation formula}\label{sec2}

Let us start with Eq. (\ref{zeta1}):
\begin{equation}
\ln A=\frac{1}{12}-\zeta'(-1).
\end{equation}
Using Eq. (\ref{derivzeta}), we get
\begin{equation}
\zeta'(-1)=-\frac{1}{4}+\frac{1}{12}+\frac{1}{6}I_3(-1)=\frac{1}{6}\left[I_3(-1)-1\right],
\end{equation}
with
\begin{equation}
I_3(-1)=\int_1^{\infty}\frac{P_3(x)}{x^{2}}\,dx,
\end{equation}
where
\begin{equation}
P_3(x)=\frac{3}{2\pi^3}\sum_{k=1}^{\infty}\frac{\sin(2k\pi x)}{k^3}.
\end{equation}
The integral $I_3(-1)$ can be put in the form
\begin{equation}
I_3(-1)=\frac{3i}{4\pi^3}\left[\int_1^{\infty}\frac{\mathrm{Li}_3(e^{-2i\pi x})}{x^2}\,dx-\int_1^{\infty}\frac{\mathrm{Li}_3(e^{2i\pi x})}{x^2}\,dx\right],
\end{equation}
involving the trilogarithm, \emph{i.e.} the polylogarithm
\begin{equation}
\mathrm{Li}_n(z)=\sum_{k=1}^{\infty}\frac{z^k}{k^n}
\end{equation}
of order 3. One obtains
\begin{equation}
\int_1^{\infty}\frac{\sin(2k\pi x)}{x^2}\,dx=-2\pi k\mathrm{Ci}(2k\pi),
\end{equation}
yielding
\begin{equation}
I_3(-1)=-\frac{3}{\pi^2}\sum_{k=1}^{\infty}\frac{\mathrm{Ci}(2k\pi)}{k^2}
\end{equation}
and thus
\begin{equation}
\zeta'(-1)=-\frac{1}{6}-\frac{1}{2\pi^2}\sum_{k=1}^{\infty}\frac{\mathrm{Ci}(2k\pi)}{k^2},
\end{equation}
where $\mathrm{Ci}$ is the usual cosine integral
\begin{equation}
\mathrm{Ci}(z)=-\int_z^{\infty}\frac{\cos t}{t}dt,    
\end{equation}
which has the following series expansion
\begin{equation}
\mathrm{Ci}(z)=\gamma+\ln z+\sum_{k=1}^{\infty}\frac{(-z^2)^k}{2k(2k)!}. 
\end{equation} 
One gets {\it in fine}
\begin{empheq}[box=\fbox]{align}
\ln A=\frac{1}{4}\left[1+\frac{2}{\pi^2}\sum_{k=1}^{\infty}\frac{\mathrm{Ci}(2k\pi)}{k^2}\right].
\end{empheq}

\section{Second expression}\label{sec3}

Let us now consider the expression (Eq. (\ref{zeta2})):
\begin{equation}
\ln A=\frac{1}{12}\ln(2\pi)+\frac{\gamma}{12}-\frac{\zeta'(2)}{2\pi^2}.
\end{equation}
In virtue of Eq. (\ref{derivzeta}), one has
\begin{equation}
\zeta'(2)=-1+\frac{1}{12}-\frac{24}{6}I'_3(2)-\frac{13}{3}I_3(2),
\end{equation}
\emph{i.e.}
\begin{equation}
\zeta'(2)=\frac{11}{12}-4I'_3(2)-\frac{13}{3}I_3(2),
\end{equation}
with
\begin{equation}
I_3(2)=\int_1^{\infty}\frac{P_3(x)}{x^5}\,dx
\end{equation}
and
\begin{equation}
I_3(2)=\frac{3i}{4\pi^3}\left[\int_1^{\infty}\frac{\mathrm{Li}_3(e^{-2i\pi x})}{x^5}\,dx-\int_1^{\infty}\frac{\mathrm{Li}_3(e^{2i\pi x})}{x^5}\,dx\right].
\end{equation}
We have the equality
\begin{equation}
\int_1^{\infty}\frac{\sin(2k\pi x)}{x^5}\,dx=\frac{1}{6}\left\{2k\pi(1-2k^2\pi^ 2)+2k^4\pi^4\left[\pi-2~\mathrm{Si}(2k\pi)\right]\right\},
\end{equation}
where $\mathrm{Si}$ represents the sine integral
\begin{equation}
\mathrm{Si}(z)=\int_0^z\frac{\sin t}{t}dt,
\end{equation}
which has the representations
\begin{equation}
\mathrm{Si}(z)=\sum_{k=0}^{\infty}\frac{(-1)^kz^{2k+1}}{(2k+1)^2(2k)!}=2x\sum_{n=0}^{\infty}\left[j_n(x)\right]^2,
\end{equation}
where $j_n$ is the spherical Bessel function of the first kind of order $n$. One has therefore 
\begin{equation}
I_3(2)=\frac{1}{4\pi^2}\sum_{k=1}^{\infty}\frac{\left\{k(1-2k^2\pi^ 2)+2k^4\pi^3(\pi-2~\mathrm{Si}(2k\pi))\right\}}{k^3}.
\end{equation}
Dealing with
\begin{equation}
I'_3(2)=-\int_1^{\infty}\frac{P_3(x)}{x^{5}}\ln x\,dx,
\end{equation}
one has also
\begin{equation}
I'_3(2)=\frac{3i}{4\pi^3}\left[\int_1^{\infty}\frac{\ln x.\mathrm{Li}_3(e^{-2i\pi x})}{x^5}\,dx-\int_1^{\infty}\frac{\mathrm{Li}_3(e^{2i\pi x})}{x^5}\,dx\right].
\end{equation}
The following equality, obtained with a Computer Algebra System \cite{Mathematica}:
\begin{eqnarray}
\int_1^{\infty}\frac{\sin(2k\pi x)}{x^5}\,dx&=&-\frac{1}{108}\left\{3k^4\pi^5[-25+12(\gamma+\ln(2\pi k)]+32k^3\pi^3~_1F_2\left[\begin{array}{c}
-\frac{1}{2}\\
\frac{1}{2}, \frac{5}{2}\\
\end{array};-k^2\pi^2\right]\right.\nonumber\\
& &\left.+24k^3\pi^3~_2F_3\left[\begin{array}{c}
-\frac{1}{2},-\frac{1}{2}\\
\frac{1}{2},\frac{1}{2},\frac{5}{2}\\
\end{array};-k^ 2\pi^2\right]\right\},
\end{eqnarray}
leads to
\begin{eqnarray}
I'_3(2)&=&\frac{1}{72}\sum_{k=1}^{\infty}\left\{3k\pi^2[-25+12(\gamma+\ln(2\pi k)]+32~_1F_2\left[\begin{array}{c}
-\frac{1}{2}\\
\frac{1}{2}, \frac{5}{2}\\
\end{array};-k^2\pi^2\right]\right.\nonumber\\
& &\left.+24~_2F_3\left[\begin{array}{c}
-\frac{1}{2},-\frac{1}{2}\\
\frac{1}{2},\frac{1}{2},\frac{5}{2}\\
\end{array};-k^2\pi^2\right]\right\},
\end{eqnarray}
and finally
\begin{empheq}[box=\fbox]{align}
\ln A&=\frac{1}{12}\left[\ln(2\pi)+\gamma+\frac{11}{2\pi^2}\right]+\sum_{k=1}^{\infty}\left[k(\gamma-1)+\frac{13}{12\pi^2}\left(\frac{1}{2\pi^2 k^2}-1\right)\right.\nonumber\\
&\left. +\frac{8}{9\pi^2}\left(~_1F_2\left[\begin{array}{c}
-\frac{1}{2}\\
\frac{1}{2}, \frac{5}{2}\\
\end{array};-k^2\pi^2\right]+3~_2F_3\left[\begin{array}{c}
-\frac{1}{2},-\frac{1}{2}\\
\frac{1}{2},\frac{1}{2},\frac{5}{2}\\
\end{array};-k^ 2\pi^2\right]\right)+k\ln(2\pi k)-\frac{13k}{6\pi}\mathrm{Si}(2k\pi)\right].
\end{empheq}

\section{Conclusion}

We derived two series expansions for the logarithm of the Glaisher-Kinkelin constant, obtained using expressions of derivatives of the Riemann zeta function. The first ones depends only on $\pi$ and cosine integrals. The second formula involves two hypergeometric functions, the Euler gamma constant, $\pi$ and sine integrals. In the future, we plan to investigate other constants of the  Glaisher-Kinkelin type \cite{Sondow2005,Choi2007,Mortici2013,Blagouchine2016,Blagouchine2018,Coppo2019,Coppo2022}, such as
\begin{equation}
\ln B=\lim_{n\rightarrow\infty}\left\{\sum_{k=1}^nk^2\ln k-\left(\frac{n^3}{3}+\frac{n^2}{2}+\frac{n}{6}\right)\ln n+\frac{n^3}{9}-\frac{n}{12}\right\}
\end{equation}
and
\begin{equation}
\ln C=\lim_{n\rightarrow\infty}\left\{\sum_{k=1}^nk^3\ln k-\left(\frac{n^4}{4}+\frac{n^3}{2}+\frac{n^2}{4}-\frac{1}{120}\right)\ln n+\frac{n^4}{16}-\frac{n^2}{12}\right\}
\end{equation}
also called sometimes the Bendersky-Adamchik constants \cite{Bendersky1933,Adamchik1998,Adamchik2005}, which were considered by Choi and Srivastava \cite{Choi1997} in the theory of multiple gamma functions.

\end{document}